\documentclass[a4paper,12pt]{article}
\usepackage[a4paper,margin=2.5truecm]{geometry}
\usepackage{amsmath,amsfonts,amssymb,amsthm,txfonts,hyperref,url,tikz}
\newtheorem{conjecture}{Conjecture}
\newtheorem{lemma}{Lemma}
\newtheorem{corollary}{Corollary}
\newtheorem{observation}{Observation}
\newtheorem{theorem}{Theorem}
\newcommand{\N}{\mathbb N}
\newcommand{\Z}{\mathbb Z}
\newcommand{\K}{\textbf{\textit{K}}}
\title{Is there a computable upper bound for the height of a solution
of a Diophantine equation with a unique solution in positive integers?}
\author{Apoloniusz Tyszka}
\begin{document}
\begin{sloppypar}
\date{}
\maketitle
\begin{abstract}
Let \mbox{$B_n=\{x_i \cdot x_j=x_k,~x_i+1=x_k:~i,j,k \in \{1,\ldots,n\}\}$}.
For a positive integer $n$, let \mbox{$\xi(n)$} denote the smallest positive
integer $b$ such that for each system \mbox{${\cal S} \subseteq B_n$}
with a unique solution in positive integers \mbox{$x_1,\ldots,x_n$}, this solution belongs to \mbox{$[1,b]^n$}.
Let \mbox{$g(1)=1$}, and let \mbox{$g(n+1)=2^{\textstyle 2^{\textstyle g(n)}}$} for every positive integer $n$.
We conjecture that \mbox{$\xi(n) \leqslant g(2n)$} for every positive integer $n$.
We prove: {\tt (1)} the function \mbox{$\xi \colon \N \setminus \{0\} \to \N \setminus \{0\}$} is
computable in the limit; {\tt (2)} if a function \mbox{$f \colon \N \setminus \{0\} \to \N \setminus \{0\}$}
has a \mbox{single-fold} Diophantine representation,
then there exists a positive integer $m$ such that \mbox{$f(n)<\xi(n)$} for every integer \mbox{$n>m$};
{\tt (3)}~the conjecture implies that there exists an algorithm which takes as input a Diophantine equation
\mbox{$D(x_1,\ldots,x_p)=0$} and returns a positive integer $d$ with the following property:
for every positive integers \mbox{$a_1,\ldots,a_p$}, if the tuple \mbox{$(a_1,\ldots,a_p)$} solely solves
the equation \mbox{$D(x_1,\ldots,x_p)=0$} in positive integers, then \mbox{$a_1,\ldots,a_p \leqslant d$};
{\tt (4)}~the conjecture implies that if a set \mbox{${\cal M} \subseteq \N$} has a \mbox{single-fold}
Diophantine representation, then ${\cal M}$ is computable; {\tt (5)}~for every integer \mbox{$n>9$},
the inequality \mbox{$\xi(n)<\left(2^{\textstyle 2^{n-5}}-1\right)^{\textstyle 2^{n-5}}+1$}
implies that \mbox{$2^{\textstyle 2^{n-5}}+1$} is composite.
\end{abstract}
\mbox{single-fold} Diophantine representation, upper bound for the height of a solution.
\section{Introduction}
In this article, we propose a conjecture which implies that there exists
an algorithm which takes as input a Diophantine equation
\mbox{$D(x_1,\ldots,x_p)=0$} and returns a positive integer $d$ with the following property:
for every positive integers \mbox{$a_1,\ldots,a_p$}, if the tuple \mbox{$(a_1,\ldots,a_p)$}
solely solves the equation \mbox{$D(x_1,\ldots,x_p)=0$} in positive integers, then \mbox{$a_1,\ldots,a_p \leqslant d$}.
Since the height of an integer tuple \mbox{$(c_1,\ldots,c_p)$} is defined as \mbox{$\max(|c_1|,\ldots,|c_p|)$},
our claim asserts that the height of the tuple \mbox{$(a_1,\ldots,a_p)$} does not exceed $d$.
\section{A conjecture on integer arithmetic}
Let \mbox{$g(1)=1$}, and let \mbox{$g(n+1)=2^{\textstyle 2^{\textstyle g(n)}}$} for every positive integer $n$
Let \mbox{$B_n=\{x_i \cdot x_j=x_k,~x_i+1=x_k:~i,j,k \in \{1,\ldots,n\}\}$}.
For a positive integer $n$, let \mbox{$\xi(n)$} denote the smallest positive
integer $b$ such that for each system \mbox{${\cal S} \subseteq B_n$} with a unique solution in
positive integers \mbox{$x_1,\ldots,x_n$}, this solution belongs to \mbox{$[1,b]^n$}.
We do not know whether or not there exists a computable function
\mbox{$\gamma \colon \N \setminus \{0\} \to \N \setminus \{0\}$} such that the inequality
\mbox{$\xi(n) \leqslant \gamma(n)$} holds for every sufficiently large positive integer $n$.
\begin{theorem}\label{the0}
The function \mbox{$\xi \colon \N \setminus \{0\} \to \N \setminus \{0\}$} is computable in the limit.
\end{theorem}
\begin{proof}
The flowchart in Figure 1 describes an algorithm which computes \mbox{$\xi(n)$} in the limit.
\vskip 0.01truecm
\begin{center}
\begin{tikzpicture}[very thick]
\tt
\node at (6,7.55) {Start};
\node at (6,6.5) {Input a positive integer $n$};
\node at (6,5.55) {$k:=1$};
\node at (6,4.1) {Create a list ${\cal C}$ of all systems ${\cal S} \subseteq B_n$ which};
\node at (6,3.6) {have exactly one solution in $\{1,\ldots,k\}^n$};
\node at (6,2.6) {Compute the smallest positive integer $b \in \{1,\ldots,k\}$ such};
\node at (6,2.1) {that each system from ${\cal C}$ has a solution in $\{1,\ldots,b\}^n$};
\node at (6,1.15) {Print $b$};
\node at (6,.25) {$k:=k+1$};

\draw (6,7.55) ellipse(.8 and .35);
\draw (2.8,6.2) -- (9,6.2) -- (9.2,6.8) -- (3,6.8) -- cycle;
\draw (5.35,5.3) rectangle (6.6,5.8);
\draw (1.2,3.3) rectangle (10.8,4.4);
\draw (0,1.8) rectangle (12,2.9);
\draw (4.9,.9) -- (6.9,.9) -- (7.1,1.4) -- (5.1,1.4) -- cycle;
\draw (5.05,0) rectangle (6.9,.5);

\draw[->] (6,7.2) -- (6,6.8);
\draw[->] (6,6.2) -- (6,5.8);
\draw[->] (6,5.3) -- (6,4.4);
\draw[->] (6,3.3) -- (6,2.9);
\draw[->] (6,1.8) -- (6,1.4);
\draw[->] (6,.9) -- (6,.5);
\draw[->] (6.9,.25) -- (12.3,.25) -- (12.3,4.9) -- (6,4.9);
\end{tikzpicture}
\end{center}
\vskip 0.01truecm
\centerline{{\bf Fig.~1}~~An infinite computation of \mbox{$\xi(n)$}}
\end{proof}
\begin{lemma}\label{lem1}
For every positive integers $b$ and $c$, \mbox{$b+1=c$} if and only if
\mbox{$2^{\textstyle 2^b} \cdot 2^{\textstyle 2^b}=2^{\textstyle 2^c}$}.
\end{lemma}
Observations~\ref{obs1} and \ref{obs2} heuristically justify Conjecture~\ref{con1}.
\begin{observation}\label{obs1}
For every system \mbox{${\cal S} \subseteq B_n$}, the following new system
\[
\left(\bigcup_{x_i \cdot x_j=x_k \in S} \{x_i \cdot x_j=x_k\}\right) \cup
\{2^{\textstyle 2^{\textstyle x_k}}=y_k:~k \in \{1,\ldots,n\}\} \cup \left(\bigcup_{x_i+1=x_k \in S} \{y_i \cdot y_i=y_k\}\right)
\]
is equivalent to ${\cal S}$.
If the system ${\cal S}$ has exactly one solution in positive integers \mbox{$x_1,\ldots,x_n$},
then the new system has exactly one solution in positive integers \mbox{$x_1,\ldots,x_n,y_1,\ldots,y_n$}.
\end{observation}
\begin{proof}
It follows from Lemma~\ref{lem1}.
\end{proof}
\begin{observation}\label{obs2}
For every positive integer $n$, the following system
\[
\left\{\begin{array}{rcl}
x_1 \cdot x_1 &=& x_1 \\
\forall i \in \{1,\ldots,n-1\} ~2^{\textstyle 2^{\textstyle x_i}} &=& x_{i+1} ~({\rm if~} n>1)
\end{array}\right.
\]
has exactly one solution in positive integers, namely \mbox{$(g(1),\ldots,g(n))$}.
\end{observation}
\par
Lemma~\ref{lem1} and \mbox{Observations~\ref{obs1}--\ref{obs2}} justify the following conjecture.
\begin{conjecture}\label{con1}
If a system \mbox{${\cal S} \subseteq B_n$} has exactly one solution in positive
integers \mbox{$x_1,\ldots,x_n$} , then \mbox{$x_1,\ldots,x_n \leqslant g(2n)$}.
In other words, \mbox{$\xi(n) \leqslant g(2n)$} for every positive integer $n$.
\end{conjecture}
\vskip 0.01truecm
\noindent
Conjecture~\ref{con2} generalizes Conjecture~\ref{con1}.
\begin{conjecture}\label{con2}
If a system \mbox{${\cal S} \subseteq B_n$} has only finitely many solutions in positive
integers \mbox{$x_1,\ldots,x_n$}, then \mbox{$x_1,\ldots,x_n \leqslant g(2n)$}.
\end{conjecture}
\section{Algebraic lemmas lead to the main theorem}
Let \mbox{$\cal{R}${\sl ng}} denote the class of all rings $\K$ that extend $\Z$, and let
\[
E_n=\{1=x_k,~x_i+x_j=x_k,~x_i \cdot x_j=x_k:~i,j,k \in \{1,\ldots,n\}\}
\]
\mbox{Th. Skolem} proved that every Diophantine equation can be algorithmically transformed
into an equivalent system of Diophantine equations of degree at most~$2$,
see \mbox{\cite[pp.~2--3]{Skolem}} and \mbox{\cite[pp.~3--4]{Matiyasevich1}}.
The following result strengthens Skolem's theorem.
\begin{lemma}\label{lem2} (\cite[p.~720]{Tyszka2})
Let \mbox{$D(x_1,\ldots,x_p) \in {\Z}[x_1,\ldots,x_p]$}.
Assume that \mbox{${\rm deg}(D,x_i) \geqslant 1$} for each \mbox{$i \in \{1,\ldots,p\}$}. We can compute a positive
integer \mbox{$n>p$} and a system \mbox{${\cal T} \subseteq E_n$} which satisfies the following two conditions:
\vskip 0.2truecm
\noindent
{\tt Condition 1.} If \mbox{$\K \in {\cal R}{\sl ng} \cup \{\N,~\N \setminus \{0\}\}$}, then
\[
\forall \tilde{x}_1,\ldots,\tilde{x}_p \in \K ~\Bigl(D(\tilde{x}_1,\ldots,\tilde{x}_p)=0 \Longleftrightarrow
\exists \tilde{x}_{p+1},\ldots,\tilde{x}_n \in \K ~(\tilde{x}_1,\ldots,\tilde{x}_p,\tilde{x}_{p+1},\ldots,\tilde{x}_n) ~solves~ {\cal T}\Bigr)
\]
{\tt Condition 2.} If \mbox{$\K \in {\cal R}{\sl ng} \cup \{\N,~\N \setminus \{0\}\}$}, then
for each \mbox{$\tilde{x}_1,\ldots,\tilde{x}_p \in \K$} with \mbox{$D(\tilde{x}_1,\ldots,\tilde{x}_p)=0$},
there exists a unique tuple \mbox{$(\tilde{x}_{p+1},\ldots,\tilde{x}_n) \in {\K}^{n-p}$} such that the tuple
\mbox{$(\tilde{x}_1,\ldots,\tilde{x}_p,\tilde{x}_{p+1},\ldots,\tilde{x}_n)$} solves ${\cal T}$.
\vskip 0.2truecm
\noindent
Conditions 1 and 2 imply that for each \mbox{$\K \in {\cal R}{\sl ng} \cup \{\N,~\N \setminus \{0\}\}$},
the equation \mbox{$D(x_1,\ldots,x_p)=0$} and the system ${\cal T}$ have the same number of solutions in $\K$.
\end{lemma}
\par
Let $\alpha$, $\beta$, and $\gamma$ denote variables.
\begin{lemma}\label{lem3} (\cite[p.~100]{Robinson})
For each positive integers \mbox{$x,y,z$}, \mbox{$x+y=z$} if and only if
\[
(zx+1)(zy+1)=z^2(xy+1)+1
\]
\end{lemma}
\begin{corollary}\label{cor1}
We can express the equation \mbox{$x+y=z$} as an equivalent system ${\cal F}$,
where ${\cal F}$ involves \mbox{$x,y,z$} and $9$ new variables, and where ${\cal F}$ consists of equations
of the forms \mbox{$\alpha+1=\gamma$} and \mbox{$\alpha \cdot \beta=\gamma$}.
\end{corollary}
\begin{proof}
The new $9$ variables express the following polynomials:
\[
zx,~~~~~zx+1,~~~~~zy,~~~~~zy+1,~~~~~z^2,~~~~~xy,~~~~~xy+1,~~~~~~z^2(xy+1),~~~~~z^2(xy+1)+1
\]
\end{proof}
\begin{lemma}\label{lem4}
Let \mbox{$D(x_1,\ldots,x_p) \in {\Z}[x_1,\ldots,x_p]$}.
Assume that \mbox{${\rm deg}(D,x_i) \geqslant 1$} for each \mbox{$i \in \{1,\ldots,p\}$}. We can compute a positive
integer \mbox{$n>p$} and a system \mbox{${\cal T} \subseteq B_n$} which satisfies the following two conditions:
\vskip 0.2truecm
\noindent
{\tt Condition 3.} For every positive integers \mbox{$\tilde{x}_1,\ldots,\tilde{x}_p$},
\[
D(\tilde{x}_1,\ldots,\tilde{x}_p)=0 \Longleftrightarrow
\exists \tilde{x}_{p+1},\ldots,\tilde{x}_n \in \N \setminus \{0\} ~(\tilde{x}_1,\ldots,\tilde{x}_p,\tilde{x}_{p+1},\ldots,\tilde{x}_n) ~solves~ {\cal T}
\]
{\tt Condition 4.} If positive integers \mbox{$\tilde{x}_1,\ldots,\tilde{x}_p$} satisfy
\mbox{$D(\tilde{x}_1,\ldots,\tilde{x}_p)=0$}, then there exists a unique tuple
\mbox{$(\tilde{x}_{p+1},\ldots,\tilde{x}_n) \in (\N \setminus \{0\})^{n-p}$} such that the tuple
\mbox{$(\tilde{x}_1,\ldots,\tilde{x}_p,\tilde{x}_{p+1},\ldots,\tilde{x}_n)$} solves ${\cal T}$.
\vskip 0.2truecm
\noindent
Conditions 3 and 4 imply that the equation \mbox{$D(x_1,\ldots,x_p)=0$} and the system ${\cal T}$ have
the same number of solutions in positive integers.
\end{lemma}
\begin{proof}
Let the system ${\cal T}$ be given by Lemma~\ref{lem2}.
We replace in ${\cal T}$ each equation of the form \mbox{$1=x_k$} by the equation \mbox{$x_k \cdot x_k=x_k$}.
Next, we apply Corollary~\ref{cor1} and replace in ${\cal T}$ each equation of the form \mbox{$x_i+x_j=x_k$}
by an equivalent system of equations of the forms \mbox{$\alpha+1=\gamma$} and \mbox{$\alpha \cdot \beta=\gamma$}.
\end{proof}
\par
Lemma~\ref{lem4} implies Theorem~\ref{the1}.
\begin{theorem}\label{the1}
Conjecture~\ref{con1} implies that there exists an algorithm which takes as input a Diophantine equation
\mbox{$D(x_1,\ldots,x_p)=0$} and returns a positive integer $d$ with the following property:
for every positive integers \mbox{$a_1,\ldots,a_p$}, if the tuple \mbox{$(a_1,\ldots,a_p)$} solely solves
the equation \mbox{$D(x_1,\ldots,x_p)=0$} in positive integers, then \mbox{$a_1,\ldots,a_p \leqslant d$}.
\end{theorem}
\begin{theorem}\label{the2}
Conjecture~\ref{con1} implies that there exists an algorithm which takes as input a Diophantine equation
\mbox{$D(x_1,\ldots,x_p)=0$} and returns a positive integer $d$ with the following property:
for every \mbox{non-negative} integers \mbox{$a_1,\ldots,a_p$}, if the tuple \mbox{$(a_1,\ldots,a_p)$} solely
solves the equation \mbox{$D(x_1,\ldots,x_p)=0$} in \mbox{non-negative} integers, then \mbox{$a_1,\ldots,a_p \leqslant d$}.
\end{theorem}
\begin {proof}
We apply Theorem~\ref{the1} to the equation \mbox{$D(x_{1}-1,\ldots,x_{p}-1)=0$}.
\end{proof}
\section{Single-fold Diophantine representations}
The Davis-Putnam-Robinson-Matiyasevich theorem states that every recursively
enumerable set \mbox{${\cal M} \subseteq {\N}^n$} has a Diophantine
representation, that is
\[
(a_1,\ldots,a_n) \in {\cal M} \Longleftrightarrow \exists x_1, \ldots, x_m \in \N ~~W(a_1,\ldots,a_n,x_1,\ldots,x_m)=0 \tag*{\texttt{(R)}}
\]
for some polynomial $W$ with integer coefficients, see \cite{Matiyasevich1}.
The polynomial~$W$ can be computed, if we know the Turing \mbox{machine $M$} such
that, for all \mbox{$(a_1,\ldots,a_n) \in {\N}^n$}, $M$ halts on \mbox{$(a_1,\ldots,a_n)$} if
and only if \mbox{$(a_1,\ldots,a_n) \in {\cal M}$}, \mbox{see \cite{Matiyasevich1}}.
The representation~\texttt{(R)} is said to be \mbox{single-fold}, if for any
\mbox{$a_1,\ldots,a_n \in \N$} the equation
\mbox{$W(a_1,\ldots,a_n,x_1,\ldots,x_m)=0$} has at most one solution
\mbox{$(x_1,\ldots,x_m) \in {\N}^m$}.
The representation~\texttt{(R)} is said to be \mbox{finite-fold}, if for any
\mbox{$a_1,\ldots,a_n \in \N$} the equation
\mbox{$W(a_1,\ldots,a_n,x_1,\ldots,x_m)=0$} has only finitely many solutions
\mbox{$(x_1,\ldots,x_m) \in {\N}^m$}. \mbox{Yu. Matiyasevich} conjectured that
each recursively enumerable set \mbox{${\cal M} \subseteq {\N}^n$} has a
\mbox{single-fold} (\mbox{finite-fold}) Diophantine representation, see \mbox{\cite[pp.~341--342]{DMR}} and
\mbox{\cite[p.~42]{Matiyasevich2}}. Currently, he seems agnostic on his conjectures,
see \mbox{\cite[p.~749]{Matiyasevich3}}.
In \mbox{\cite[p.~581]{Tyszka1}}, the author explains why Matiyasevich's conjectures although widely known are less widely accepted.
\begin{theorem}\label{the3}~(cf.~\cite[Theorem~5,~p.~711]{Tyszka0})
Conjecture~\ref{con1} implies that if a set \mbox{${\cal M} \subseteq \N$} has a \mbox{single-fold}
Diophantine representation, then ${\cal M}$ is computable.
In particular, Conjecture~\ref{con1} contradicts Matiyasevich's
conjecture on \mbox{single-fold} Diophantine representations.
\end{theorem}
\begin{proof}
Let a set \mbox{${\cal M} \subseteq \N$} have a \mbox{single-fold} Diophantine representation.
It means that there exists a polynomial \mbox{$W(x,x_1,\ldots,x_m)$} with integer coefficients such that
\[
\forall b \in \N~ \Bigl(b \in {\cal M} \Longleftrightarrow \exists x_1, \ldots, x_m \in \N ~~W(b,x_1,\ldots,x_m)=0\Bigr)
\]
and for every \mbox{$b \in \N$} the equation \mbox{$W(b,x_1,\ldots,x_m)=0$} has at most one solution
\mbox{$(x_1,\ldots,x_m) \in {\N}^m$}. By Theorem~\ref{the2}, there exists a computable function \mbox{$\theta \colon \N \to \N$}
such that for every \mbox{$b,x_1,\ldots,x_m \in \N$} the equality \mbox{$W(b,x_1,\ldots,x_m)=0$} implies
\mbox{$\max(x_1,\ldots,x_m) \leqslant \theta(b)$}.
Hence, we can decide whether or not a \mbox{non-negative} integer $b$ belongs to ${\cal M}$ by checking
whether or not the equation \mbox{$W(b,x_1,\ldots,x_m)=0$} has an integer solution in the box \mbox{$[0,\theta(b)]^m$}.
\end{proof}
\begin{observation}\label{obs3}
Theorem~\ref{the3} remains true if we change the bound \mbox{$g(2n)$} in Conjecture~\ref{con1}
to any other computable bound \mbox{$\delta(n)$}.
\end{observation}
\begin{theorem}\label{the4}
If a function \mbox{$f \colon \N \setminus \{0\} \to \N \setminus \{0\}$}
has a \mbox{single-fold} Diophantine representation,
then there exists a positive integer $m$ such that \mbox{$f(n)<\xi(n)$} for every integer \mbox{$n>m$}.
\end{theorem}
\begin{proof}
There exists a polynomial \mbox{$W(x_1,x_2,x_3,\ldots,x_r)$} with integer coefficients
such that for each positive integers \mbox{$x_1,x_2$},
\[
(x_1,x_2) \in f \Longleftrightarrow \exists x_3,\ldots,x_r \in \N \setminus \{0\}~~ W(x_1,x_2,x_3-1,\ldots,x_r-1)=0
\]
and for each positive integers \mbox{$x_1,x_2$} at most one tuple
\mbox{$(x_3,\ldots,x_r)$} of positive integers satisfies \mbox{$W(x_1,x_2,x_3-1,\ldots,x_r-1)=0$}.
By Lemma~\ref{lem4}, there exists an integer \mbox{$s \geqslant 3$} such that
for every positive integers \mbox{$x_1,x_2$},
\begin{equation}
(x_1,x_2) \in f \Longleftrightarrow \exists x_3,\ldots,x_s \in \N \setminus \{0\}~~ \Psi(x_1,x_2,x_3,\ldots,x_s)\tag*{\texttt{(E)}}
\end{equation}
where \mbox{$\Psi(x_1,x_2,x_3,\ldots,x_s)$} is a conjunction of formulae of the forms
\mbox{$x_i+1=x_k$} and \mbox{$x_i \cdot x_j=x_k$}, the indices $i,j,k$ belong to
$\{1,\ldots,s\}$, and for each positive integers \mbox{$x_1,x_2$} at most one
tuple \mbox{$(x_3,\ldots,x_s)$} of positive integers satisfies \mbox{$\Psi(x_1,x_2,x_3,\ldots,x_s)$}.
Let $[\cdot]$ denote the integer part function, and let an integer~$n$ is greater than \mbox{$m=2s+2$}.
Then,
\[
n \geqslant \left[\frac{n}{2}\right]+\frac{n}{2}>\left[\frac{n}{2}\right]+s+1
\]
and \mbox{$n-\left[\frac{n}{2}\right]-s-2 \geqslant 0$}.
Let $T_n$ denote the following system with $n$ variables:
\[
\left\{
\begin{array}{c}
\textrm{all~equations~occurring~in~}\Psi(x_1,x_2,x_3,\ldots,x_s)\\
\begin{array}{rcl}
\forall i \in \left\{1,\ldots,n-\left[\frac{n}{2}\right]-s-2\right\} ~u_i \cdot u_i &=& u_i\\
t_1 \cdot t_1 &=& t_1\\
\forall i \in \left\{1,\ldots,\left[\frac{n}{2}\right]-1\right\} ~t_i+1 &=& t_{i+1}\\
t_2 \cdot t_{\left[\frac{n}{2}\right]} &=& u\\
u+1 &=&x_1 {\rm ~(if~}n{\rm ~is~odd)}\\
t_1 \cdot u &=& x_1 {\rm ~(if~}n{\rm ~is~even)}\\
x_2+1 &=& y
\end{array}
\end{array}
\right.
\]
By the equivalence~\texttt{(E)}, the system \mbox{$T_n$} is solvable in positive integers,
\mbox{$2 \cdot \left[\frac{n}{2}\right]=u$}, \mbox{$n=x_1$}, and
\[
f(n)=f(x_1)=x_2<x_2+1=y
\]
Since \mbox{$T_n \subseteq B_n$} and \mbox{$T_n$} has at most one solution in positive
integers, \mbox{$y \leqslant \xi(n)$}. Hence, \mbox{$f(n)<\xi(n)$}.
\end{proof}
\begin{corollary}\label{cor2}
If the function \mbox{$\xi \colon \N \setminus \{0\} \to \N \setminus \{0\}$}
is dominated by a computable function, then Matiyasevich's conjecture on single-fold Diophantine
representations is false. In particular, Conjecture~\ref{con1} contradicts Matiyasevich's
conjecture on \mbox{single-fold} Diophantine representations.
\end{corollary}
\section{Fermat primes}
\begin{observation}\label{obs4}
Only \mbox{$x_1=1$} solves the equation \mbox{$x_1 \cdot x_1=x_1$} in positive integers.
Only \mbox{$x_1=1$} and \mbox{$x_2=2$} solve the system \mbox{$\{x_1 \cdot x_1=x_1,~x_1+1=x_2\}$}
in positive integers. For each integer \mbox{$n \geqslant 3$}, the following system
\begin{displaymath}
\left\{
\begin{array}{rcl}
x_1 \cdot x_1 &=& x_1 \\
x_1+1 &=& x_2 \\
\forall i \in \{2,\ldots,n-1\} ~x_i \cdot x_i &=& x_{i+1}
\end{array}
\right.
\end{displaymath}
\noindent
has a unique solution in positive integers,
namely \mbox{$\left(1,2,4,16,256,\ldots,2^{\textstyle 2^{n-3}},2^{\textstyle 2^{n-2}}\right)$}.
\end{observation}
\begin{corollary}\label{cor3}
We have: \mbox{$\xi(1)=1$} and \mbox{$\xi(2)=2$}.
The inequality \mbox{$\xi(n) \geqslant 2^{\textstyle 2^{n-2}}$}
holds for every integer \mbox{$n \geqslant 3$}.
\end{corollary}
\par
Primes of the form \mbox{$2^{\textstyle 2^n}+1$} are called Fermat primes,
as Fermat conjectured that every integer of the form \mbox{$2^{\textstyle 2^n}+1$} is prime (\mbox{\cite[p.~1]{17lectures}}).
Fermat correctly remarked that \mbox{$2^{\textstyle 2^0}+1=3$}, \mbox{$2^{\textstyle 2^1}+1=5$},
\mbox{$2^{\textstyle 2^2}+1=17$}, \mbox{$2^{\textstyle 2^3}+1=257$}, and \mbox{$2^{\textstyle 2^4}+1=65537$}
are all prime (\mbox{\cite[p.~1]{17lectures}}). It is not known whether or not there are any other Fermat
primes (\mbox{\cite[p.~206]{17lectures}}).
\begin{theorem}\label{the5}
If $n \in {\mathbb N}\setminus\{0\}$ and $2^{\textstyle 2^n}+1$ is prime,
then the following system
\[
\left\{\begin{array}{rcl}
\forall i \in \{1,\ldots,n\} ~x_i \cdot x_i &=& x_{i+1} \\
x_{1}+1 &=& x_{n+2} \\
x_{n+2}+1 &=& x_{n+3} \\
x_{n+1}+1 &=& x_{n+4} \\
x_{n+3} \cdot x_{n+5} &=& x_{n+4}
\end{array}\right.
\]
has a unique solution $\left(a_1,\ldots,a_{n+5}\right)$ in non-negative
integers. The numbers $a_1,\ldots,a_{n+5}$ are positive and
${\rm max}\left(a_1,\ldots,a_{n+5}\right)=a_{n+4}=\left(2^{\textstyle 2^n}-1\right)^{\textstyle 2^n}+1$.
\end{theorem}
\begin{proof}
The system equivalently expresses that \mbox{$(x_1+2) \cdot x_{n+5}=x_1^{\textstyle 2^n}+1$}.
Since
\[
x_1^{\textstyle 2^n}+1=\left((x_1+2)-2\right)^{\textstyle 2^n}+1=
\]
\[
2^{\textstyle 2^n}+1+(x_1+2) \cdot \sum_{\textstyle k=1}^{\textstyle 2^n}
{\textstyle 2^n \choose k} \cdot (x_1+2)^{\textstyle k-1} \cdot (-2)^{\textstyle 2^n-k}
\]
we get
\[
(x_1+2) \cdot \left(x_{n+5}-\sum_{\textstyle k=1}^{\textstyle 2^n} {{\textstyle 2^n} \choose {\textstyle k}}
\cdot (x_1+2)^{\textstyle k-1} \cdot (-2)^{\textstyle 2^n-k}\right)=2^{\textstyle 2^n}+1
\]
Hence, $x_1+2$ divides $2^{\textstyle 2^n}+1$.
Since $x_1+2 \geqslant 2$ and $2^{\textstyle 2^n}+1$ is prime, we get
$x_1+2=2^{\textstyle 2^n}+1$ and $x_1=2^{\textstyle 2^n}-1$. Next,
$x_{n+1}=x_1^{\textstyle 2^n}=\left(2^{\textstyle 2^n}-1\right)^{\textstyle 2^n}$ and
\[
x_{n+4}=x_{n+1}+1=\left(2^{\textstyle 2^n}-1\right)^{\textstyle 2^n}+1
\]
Explicitly, the whole solution is given by
\[
\left\{\begin{array}{rcl}
\forall i \in \{1,\ldots,n+1\} ~a_i &=& \left(2^{\textstyle 2^n}-1\right)^{\textstyle 2^{i-1}} \\
a_{n+2} &=& 2^{\textstyle 2^n} \\
a_{n+3} &=& 2^{\textstyle 2^n}+1 \\
a_{n+4} &=& \left(2^{\textstyle 2^n}-1\right)^{\textstyle 2^n}+1 \\
a_{n+5} &=& 1+\sum_{\textstyle k=1}^{\textstyle 2^{n}}\limits
\displaystyle {{\textstyle 2^n} \choose {\textstyle k}} \cdot \left(2^{\textstyle 2^n}+1\right)^{\textstyle k-1}
\cdot (-2)^{\textstyle 2^n-k}
\end{array}\right.
\]
\end{proof}
\begin{corollary}\label{cor4}
For every integer \mbox{$n>5$}, if \mbox{$2^{\textstyle 2^{n-5}}+1$} is prime, then
\[
\xi(n) \geqslant \left(2^{\textstyle 2^{n-5}}-1\right)^{\textstyle 2^{n-5}}+1
\]
In particular, 
\[
\xi(9) \geqslant \left(2^{\textstyle 2^{9-5}}-1\right)^{\textstyle 2^{9-5}}+1=\left(2^{16}-1\right)^{16}+1>
\left(2^{16}-2^{15}\right)^{16}=\left(2^{15}\right)^{16}=2^{240}>2^{\textstyle 2^{9-2}}
\]
The numbers \mbox{$2^{\textstyle 2^{n-5}}+1$} are prime when \mbox{$n \in \{6,7,8\}$}, but
\[
\left(2^{\textstyle 2^{6-5}}-1\right)^{\textstyle 2^{6-5}}+1=10<65536=2^{\textstyle 2^{6-2}}
\]
\[
\left(2^{\textstyle 2^{7-5}}-1\right)^{\textstyle 2^{7-5}}+1=50626<4294967296=2^{\textstyle 2^{7-2}}
\]
\[
\left(2^{\textstyle 2^{8-5}}-1\right)^{\textstyle 2^{8-5}}+1=17878103347812890626<18446744073709551616=2^{\textstyle 2^{8-2}}
\]
\end{corollary}
\begin{corollary}\label{cor5}
For every integer \mbox{$n>9$}, the inequality \mbox{$\xi(n)<\left(2^{\textstyle 2^{n-5}}-1\right)^{\textstyle 2^{n-5}}+1$}
implies that \mbox{$2^{\textstyle 2^{n-5}}+1$} is composite.
\end{corollary}
\section{Two stronger conjectures}
Conjecture~\ref{con3} strengthens Conjecture~\ref{con1}.
\begin{conjecture}\label{con3}
If \mbox{$n \in \{3,\ldots,8\}$}, then \mbox{$\xi(n)=2^{\textstyle 2^{n-2}}$}.
If \mbox{$n \geqslant 9$}, then
\mbox{$\xi(n) \leqslant \left(2^{\textstyle 2^{n-5}}-1\right)^{\textstyle 2^{n-5}}+1$}.
\end{conjecture}
\begin{theorem}\label{the6} (cf.~\cite[Theorem~2,~p.~709]{Tyszka0})
For each positive integer $n$, the following system
\[
\left\{\begin{array}{rcl}
\forall i \in \{1,\ldots,n\} ~x_i \cdot x_i &=& x_{i+1} \\
x_{n+2}+1 &=& x_1 \\
x_{n+3}+1 &=& x_{n+2} \\
x_{n+3} \cdot x_{n+4} &=& x_{n+1} \\
x_{n+5} \cdot x_{n+5} &=& x_{n+5} \\
x_{n+5}+1 &=& x_{n+6} \\
x_{n+6} \cdot x_{n+7} &=& x_1 \\
x_{n+6} \cdot x_{n+8} &=& x_{n+9} \\
x_{n+9}+1 &=& x_{n+4}
\end{array}\right.
\]
has a unique solution \mbox{$(a_1,\ldots,a_{n+9})$} in positive integers.
The numbers \mbox{$a_1,\ldots,a_{n+9}$} satisfy:
\[
\begin{array}{rcl}
\forall i \in \{1,\ldots,n+1\} ~a_i &=& \left(2+2^{\textstyle 2^n}\right)^{\textstyle 2^{i-1}} \\
a_{n+2} &=& 1+2^{\textstyle 2^n} \\
a_{n+3} &=& 2^{\textstyle 2^n} \\
a_{n+4} &=& \left(1+2^{\textstyle 2^n-1}\right)^{\textstyle 2^n} \\
a_{n+5} &=& 1 \\
a_{n+6} &=& 2 \\
a_{n+7} &=& 1+2^{\textstyle 2^n-1} \\
a_{n+8} &=& \frac{\textstyle \left(1+2^{\textstyle 2^n-1}\right)^{\textstyle 2^n}-1}{\textstyle 2} \\
a_{n+9} &=& \left(1+2^{\textstyle 2^n-1}\right)^{\textstyle 2^n}-1
\end{array}
\]
\end{theorem}
\begin{proof}
The tuple \mbox{$(a_1,\ldots,a_{n+9})$} consists of positive integers and solves the system.
We prove that this solution is unique in positive integers. The equations
\[
\begin{array}{rcl}
x_{n+5} \cdot x_{n+5} &=& x_{n+5} \\
x_{n+5}+1 &=& x_{n+6} \\
x_{n+6} \cdot x_{n+7} &=& x_1
\end{array}
\]
imply that $x_1$ is even. The equations
\[
\begin{array}{rcl}
x_{n+5} \cdot x_{n+5} &=& x_{n+5} \\
x_{n+5}+1 &=& x_{n+6} \\
x_{n+6} \cdot x_{n+8} &=& x_{n+9} \\
x_{n+9}+1 &=& x_{n+4}
\end{array}
\]
imply that $x_{n+4}$ is odd. The equations
\[
\begin{array}{rcl}
\forall i \in \{1,\ldots,n\} ~x_i \cdot x_i &=& x_{i+1} \\
x_{n+2}+1 &=& x_1 \\
x_{n+3}+1 &=& x_{n+2} \\
x_{n+3} \cdot x_{n+4} &=& x_{n+1}
\end{array}
\]
imply that \mbox{$x_1=x_{n+3}+2 \geqslant 3$} and
\mbox{$x_1^{\textstyle 2^n}=(x_1-2) \cdot x_{n+4}$}. This equality and the polynomial identity
\[
x_1^{\textstyle 2^n}=2^{\textstyle 2^n}+(x_1-2) \cdot
\left(2^{\textstyle 2^n-1}+\sum_{\textstyle k=1}^{\textstyle 2^n-1} 2^{\textstyle 2^n-1-k} \cdot x_1^k\right)
\]
imply that \mbox{$x_1-2$} divides $2^{\textstyle 2^n}$ and
\begin{equation}
x_{n+4}=\frac{\textstyle x_1^{\textstyle 2^n}}{\textstyle x_1-2}=
\frac{\textstyle 2^{\textstyle 2^n}}{\textstyle x_1-2}+\underbrace{2^{\textstyle 2^n-1}+\sum_{\textstyle k=1}^{\textstyle 2^n-1}
2^{\textstyle 2^n-1-k} \cdot x_1^k}_{\textstyle {{\rm the~sum~of~two~even~numbers~as~}x_1{\rm ~is~even}}}\tag*{\texttt{(L)}}
\end{equation}
Since \mbox{$x_1 \geqslant 3$} and
\mbox{$x_1-2$} divides \mbox{$2^{\textstyle 2^n}$}, \mbox{$x_1=2+2^k$}, where \mbox{$k \in \left[0,2^n\right] \cap \Z$}.
Since $x_{n+4}$ is odd, the equality $\texttt{(L)}$ implies that \mbox{$\frac{\textstyle 2^{\textstyle 2^n}}{\textstyle x_1-2}$} is odd.
Hence, \mbox{$x_1=2+2^{\textstyle 2^n}$}. Consequently, \mbox{$x_i=a_i$} for every \mbox{$i \in \{1,\ldots,n+9\}$}.
\end{proof}
\begin{corollary}\label{cor6}
The inequality \mbox{$\xi(n) \geqslant \left(2+2^{\textstyle 2^{n-9}}\right)^{\textstyle 2^{n-9}}$}
holds for every integer \mbox{$n \geqslant 10$}.
\end{corollary}
\begin{conjecture}\label{con4}
The equality \mbox{$\xi(n)=\left(2+2^{\textstyle 2^{n-9}}\right)^{\textstyle 2^{n-9}}$} holds
for every sufficiently large positive integer $n$.
\end{conjecture}

\noindent
Apoloniusz Tyszka\\
University of Agriculture\\
Faculty of Production and Power Engineering\\
Balicka 116B, 30-149 Krak\'ow, Poland\\
Email: \url{rttyszka@cyf-kr.edu.pl}
\end{sloppypar}
\end{document}